\documentclass[a4paper,11pt]{amsart}
\usepackage[T1]{fontenc}
\usepackage[latin1]{inputenc}
\usepackage{amssymb, amsmath, amsfonts, amsthm}
\usepackage{url}
\usepackage{euscript}
\usepackage[all]{xy}

\DeclareMathOperator{\id}{Id}

\DeclareMathOperator{\supp}{supp}
\newcommand{\dott}{\, \cdot\,}
\newcommand{\epsi}{\varepsilon}

\newcommand{\quot}{{\F/\Gr}}
\newcommand{\Gr}{G}

\newcommand{\N}{\ensuremath{\mathcal{N}}}

\newcommand{\D}{\ensuremath{\mathcal{D}}}
\newcommand{\G}{\ensuremath{\mathcal{G}}}
\newcommand{\F}{\ensuremath{\mathcal{F}}}
\newcommand{\C}{\ensuremath{\mathcal{C}}}
\newcommand{\R}{\ensuremath{\mathcal{R}}}
\newcommand{\tS}{\ensuremath{\tilde S}}

\newcommand{\abs}[1]{\left\vert#1\right\vert}
\newcommand{\scal}[1]{\left<#1\right>}
\newcommand{\Real}{\mathbb R}

\newcommand{\norm}[1]{\left\Vert#1\right\Vert}

\newcommand{\tnorm}[1]{\vert\hspace*{-1pt}\vert\hspace*{-1pt}\vert#1\vert\hspace*{-1pt}\vert\hspace*{-1pt}\vert}
\newcommand{\btnorm}[1]{\big|\hspace*{-1pt}\big|\hspace*{-1pt}\big|#1\big|\hspace*{-1pt}\big|\hspace*{-1pt}\big|}

\newcommand{\muac}{\mu_{\text{\rm ac}}}

\newcommand{\inv}{{^{-1}}}
\DeclareMathOperator{\sgn}{sgn}

\newtheorem{theorem}{Theorem}[section]

\newtheorem{lemma}[theorem]{Lemma}
\newtheorem{definition}[theorem]{Definition}
\newtheorem{proposition}[theorem]{Proposition}

\allowdisplaybreaks

\title{Lipschitz metric for the Hunter--Saxton equation}
\author[Bressan]{Alberto Bressan}
\address[Bressan]{Department of Mathematics, Pennsylvania State
  University, University Park, PA 16802, USA}
\email[]{bressan@math.psu.edu}
\urladdr{www.math.psu.edu/bressan/}

\author[Holden]{Helge Holden}
\address[Holden]{\newline
  Department of Mathematical Sciences,
  Norwegian University of Science and Technology,
  NO--7491 Trondheim, Norway,
  and \newline
  Centre of Mathematics for Applications,
  University of Oslo,
  P.O.\ Box 1053, Blindern,
  NO--0316 Oslo, Norway }
\email[]{holden@math.ntnu.no}
\urladdr{www.math.ntnu.no/\~{}holden/}

\author[Raynaud]{Xavier Raynaud}
\address[Raynaud]{\newline
  Department of Mathematical Sciences,
  Norwegian University of Science and Technology,
  NO--7491 Trondheim, Norway,
  and \newline
  Centre of Mathematics for Applications,
  University of Oslo,
  P.O.\ Box 1053, Blindern,
  NO--0316 Oslo, Norway}
\email[]{raynaud@math.ntnu.no}
\urladdr{www.math.ntnu.no/\~{}raynaud/}

\date{\today}
\thanks{This paper was written as part of the international research
  program on Nonlinear Partial Differential Equations at the Centre
  for Advanced Study at the Norwegian Academy of Science and Letters
  in Oslo during the year 2008--09.  Research supported in part by the Research Council of Norway.}

\begin{document}
\begin{abstract}
  We study stability of solutions of the Cauchy
  problem for the Hunter--Saxton equation
  $u_t+uu_x=\frac14(\int_{-\infty}^xu_x^2\,dx-\int_{x}^\infty
  u_x^2\,dx)$ with initial data $u_0$.  In
  particular, we derive a new Lipschitz metric $d_\D$
  with the property that for two solutions $u$ and
  $v$ of the equation we have $d_\D(u(t),v(t))\le e^{Ct}
  d_\D(u_0,v_0)$.
\end{abstract}
\maketitle
\tableofcontents

\section{Introduction}
The initial value problem for the Hunter--Saxton
equation
\begin{equation}
  \label{eq:hs0}
  u_t+uu_x=\frac14\Big(\int_{-\infty}^xu_x^2\,dx-\int_{x}^\infty
  u_x^2\,dx\Big), \quad u|_{t=0}=u_0,
\end{equation}
or alternatively
\begin{equation}
  (u_t+uu_x)_x=\frac12 u_x^2, \quad u|_{t=0}=u_0,
\end{equation}
has been widely studied since it was introduced
\cite{MR1135995} as a model for liquid crystals.
It possesses a number of startling properties,
being completely integrable, having infinitely
many conserved quantities and a Lax
pair. Furthermore, it is bi-variational and
bi-Hamiltonian \cite{MR1306466}. The initial value
problem has been extensively studied
\cite{MR1361013,MR1361014,MR1668954,MR1701136,MR1799274,MR2191785}. A
convergent finite difference scheme exists for the
equation \cite{HolKarRis:sub05}.  The simplest
conservation law reads
\begin{equation}
  (u^2_x)_t+(u u_x^2)_x=0.
\end{equation}
Furthermore, the equation
enjoys wave breaking in finite time. More precisely, if the initial
data
is {\it not} monotone
increasing, then
\begin{equation}
  \inf(u_x)\to-\infty \text{  as  } t\uparrow t^*=2/\sup(-u_0^\prime).
\end{equation}
Past wave breaking there are at least two different classes of solutions.
Consider the example \cite{MR1361014} with initial data
$u_0(x)=-x\chi_{[0,1]}(x)-\chi_{[1,\infty)}(x)$. For $t\in[0,2)$ the solution reads
\begin{equation}
  u(t,x)= \frac{2x}{t-2}\chi_{(0,(2-t)^2/4)}(x)+\frac12(t-2)\chi_{((2-t)^2/4,\infty)}(x),
  \quad t<2.
\end{equation}
Observe that $u(t,x)\to 0$ pointwise almost
everywhere as $t\to 2-$.  A careful analysis of
the solution reveals that the energy density
$u_x^2 dx$ approaches a Dirac delta mass at the
origin as $t\to 2$.  Two continuations past $t=2$
are possible: The dissipative solution
\begin{equation}
  u(t,x)=0, \quad x\in\Real,\,  t>2,
\end{equation}
and the conservative solution
\begin{equation}
  u(t,x)= \frac{2x}{t-2}\chi_{(0,(2-t)^2/4)}(x)+\frac12(t-2)\chi_{((2-t)^2/4,\infty)}(x),
  \quad t>2.
\end{equation}

Another example \cite{MR1135995}  is the following with initial data $u_0(x)=0$.
One solution (the dissipative) is clearly $u(t,x)=0$ everywhere.
Another solution (the conservative) solution reads
\begin{equation}
  u(t,x)= -2t\chi_{(-\infty,-t^2)}(x)+\frac{2x}{t}\chi_{(-t^2,t^2)}(x)+2t\chi_{(t^2,\infty)}(x).
\end{equation}
As a consequence of this the existence theory for
the Hunter--Saxton equation is complicated, and
there is a dichotomy between the dissipative and
the conservative solutions.

Zhang and Zheng \cite{MR1799274} have proved
global existence and uniqueness of both
conservative and dissipative solutions (on the
half-line $x>0$) using Young measures and
mollification techniques for compactly supported
square integrable initial data.  An alternative
approach was developed in \cite{MR2191785} for the
Hunter--Saxton equation and in \cite{BZZ} for a
somewhat more general class of nonlocal wave
equations, by rewriting the equation in terms of
an ``energy variable'', and showing the existence
of a continuous semigroup of solutions.
Furthermore, the papers \cite{MR2191785} and
\cite{BF} introduce a new distance function which
renders Lipschitz continuous this semigroup of
solutions.  This is important because it
establishes the uniqueness and continuous
dependence for the Cauchy problem.

We remark that this is a nontrivial issue for
nonlinear partial differential equations.  For
scalar conservation laws, where $u=u(t,x)\in\Real$
satisfies $u_t+\nabla_x\cdot f(u)=0$, as proved in
\cite{K} every couple of entropy weak solutions
satisfies $\|u(t)-v(t)\|_{{\bf L}^1}\le
\|u(0)-v(0)\|_{{\bf L}^1}$ for all $t\geq 0$.  For
a hyperbolic system of conservation laws in one
space dimension $u_t+f(u)_x=0$, it is well known
that, for initial data with sufficiently small
total variation, one has $\|u(t)-v(t)\|_{{\bf
    L}^1}\leq C \|u(0)-v(0)\|_{{\bf L}^1}$ for a
suitable constant $C$ and all $t$ positive
\cite{B,HR}.

The problem at hand can nicely be illustrated in
the simpler context of an ordinary differential
equation.  Consider three differential equations:
\begin{subequations}
  \begin{align}
    \dot x&= a(x), & x(0)&=x_0, \quad \text{$a$ Lipschitz},  \label{eq:reg}\\
    \dot x&= 1+\alpha H(x),&  x(0)&=x_0, \quad \text{$H$ the Heaviside function,
      $\alpha>0$},
    \label{eq:Lip}\\
    \dot x&=\abs{x}^{1/2}, & x(0)&=x_0, \quad\quad
    \text{$t\mapsto x(t)$ strictly increasing}. \label{eq:sqrtlaw}
  \end{align}
\end{subequations}
Straightforward computations give as solutions
\begin{subequations}
  \begin{align}
    x(t)&=x_0+\int_0^t a(x(s))\, ds,  \\
    x(t)&=(1+\alpha H(t-t_0))(t-t_0), \quad t_0=-x_0/(1+\alpha H(x_0)), \\
    x(t)&=\sgn\big(\frac{t}2+v_0\big) \big(\frac{t}{2}+v_0\big)^2\text{ where }v_0
    =\sgn(x_0)\abs{x_0}^{1/2}.
  \end{align}
\end{subequations}
We find that
\begin{subequations}
  \begin{align}
    \abs{x(t)-\bar x(t)}&\le e^{Lt}\abs{x_0-\bar x_0}, \quad L= \norm{a}_{\rm Lip},  \\
    \abs{x(t)-\bar x(t)}&\le (1+ \alpha) \abs{x_0-\bar x_0}, \\
    x(t)-\bar x(t)&= t v_0+\abs{x_0}, \quad \text{when $\bar x_0=0$.}
  \end{align}
\end{subequations}
Thus we see that in the regular case
\eqref{eq:reg} we get a Lipschitz estimate with
constant $e^{Lt}$ uniformly bounded as $t$ ranges
on a bounded interval.  In the second case
\eqref{eq:Lip} we get a Lipschitz estimate
uniformly valid for all $t\in \Real$.  In the
final example \eqref{eq:sqrtlaw}, by restricting
attention to strictly increasing solutions of the
ordinary differential equations, we achieve uniqueness and continuous
dependence on the initial data, but without any
Lipschitz estimate at all.  We observe that, by
introducing the Riemannian metric
\begin{equation}
  \label{eq:defriemmet}
  d(x, \bar x)= \abs{\int_x^{\bar x} \frac{dz}{{\abs{z}^{1/2}}}},
\end{equation}
an easy computation reveals that
\begin{equation}
  d(x(t), \bar x(t))= d(x_0, \bar x_0).
\end{equation}
Let us explain why this metric can be considered
as a Riemannian metric.
The Euclidean metric between the two points is
then given
\begin{equation}
  \label{eq:eucldist}
  \abs{x_0-\bar x_0}=\inf_{x}\int_0^1\abs{x_s(s)}\,ds
\end{equation}
where the infimum is taken over all paths
$x\colon[0,1]\to\Real$ that join the two points $x_0$
and $\bar x_0$, that is, $x(0)=x_0$ and
$x(1)=x_1$. However, as we have seen, the
solutions are not Lipschitz for the Euclidean
metric. Thus we want to measure the infinitesimal
variation $x_s$ in an alternative way, which makes
solutions of equation \eqref{eq:sqrtlaw} Lipschitz
continuous. We look at the evolution equation that
governs $x_s$ and, by differentiating
\eqref{eq:sqrtlaw} with respect to $s$, we get
\begin{equation*}
  \dot x_s=\frac{\sgn(x)x_s}{2\sqrt{\abs{x}}},
\end{equation*}
and we can check that
\begin{equation}
  \label{eq:dtriemevol}
  \frac{d}{dt}\left(\frac{\abs{x_s}}{\sqrt{\abs{x}}}\right)=0.
\end{equation}
Let us consider the real line as a Riemannian
manifold where, at any point $x\in\Real$, the
Riemannian norm, for any tangent vector
$v\in\Real$ in the tangent space of $x$, is given
by $\abs{v}/\sqrt{\abs{x}}$.  From
\eqref{eq:dtriemevol}, one can see that at the
infinitesimal level, this Riemannian norm is
exactly preserved by the evolution equation. The
distance on the real line which is naturally
inherited by this Riemannian is given by
\begin{equation*}
  d(x_0,\bar x_0)=\inf_{x}\int_0^1\frac{\abs{x_s}}{\sqrt{\abs{x}}}\,ds
\end{equation*}
where the infimum is taken over all paths
$x\colon[0,1]\to\Real$ joining $x_0$ and $\bar
x_0$. It is quite reasonable to restrict ourselves
to paths that satisfy $x_s\geq0$ and then, by a
change of variables, we recover the definition
\eqref{eq:defriemmet}.

We remark that, for a wide class of ordinary differential equations of the
form $\dot x = f(t,x)$, $~x\in\Real^n$, a
Riemannian metric that is contractive with respect to the
corresponding flow has been constructed in
\cite{BC90}. Here the coefficient of the metric at
a point $P=(t,x)$ depends on the total directional
variation of the (possibly discontinuous) vector
field $f$ up to the point $P$.  The equations
\eqref{eq:reg} and \eqref{eq:Lip} provide two
examples covered by this approach.

The aim of this paper is to construct a Riemannian
metric on a functional space, which renders
Lipschitz continuous the flow generated by the
Hunter--Saxton equation in the conservative
case. Let us describe the result of the paper in a
non-technical manner.  From the examples above, it
is clear that the solution itself is insufficient
to describe a unique solution. Similar to the
treatment of the Camassa--Holm equation
\cite{HolRay:07,MR2278406}, it turns out that the
appropriate way to resolve this issue to consider
the pair $(u,\mu)$ where we have added the energy
measure $\mu$ with absolute continuous part
satisfying $\muac= u_x^2 dx$. To obtain a
Lipschitz metric we introduce new variables.  To
that end assume first that one has a solution
$u=u(t,x)$, and consider the characteristics
$y_t(t,\xi)=u(t,y(t,\xi))$, the Lagrangian
velocity $U(t,\xi)=u(t,y(t,\xi))$, and the
Lagrangian cumulative energy
$H(t,\xi)=\int_{-\infty}^{y(t,\xi)}u_x^2(t,x)\,dx$.
Formally, the Hunter--Saxton equation is
equivalent to the linear system of ordinary
differential equations
\begin{equation}
  \begin{aligned} \label{eq:ode}
    y_t&=U,\\
    U_t&=\frac12H-\frac14H(\infty),\\
    H_t&=0
  \end{aligned}
\end{equation}
in an Hilbert space. The quantity
$H(\infty)=\int_{\Real}u_x^2(t,x)\,dx$ is a
constant. We first prove the existence of a global
solution, see Theorem \ref{th:semigroupS}, and the
existence of a continuous semigroup. However, in
order to return to Eulerian variables it is
necessary to resolve the redundancy, denoted
relabeling, in Lagrangian coordinates, see Section
\ref{subsec:label}. We introduce an equivalence
relation for the Lagrangian variables
corresponding to one and the same Eulerian
solution. Next, we introduce a Riemannian metric
$d$ in Lagrangian variables. Denote by
$X=(y,U,H)$. The natural choice of letting the
distance between two elements $X_0$ and $X_1$ as
the infimum of $\norm{X_0\circ f-X_1\circ f}$ over
all relabelings $f$, fails as it does not satisfy
the triangle inequality. At each point $X$, we
consider the elements that coincide to $X$ under
relabelings. Formally it corresponds to a
Riemaniann submanifold whose structure is
inherited from the ambiant Hilbert space. At each
point $X$, we show that the tangent space to the
relabeling submanifold corresponds to the set of
all elements $V$ such that $V=gX_\xi$ for some
scalar function $g$. Given $X$ and a tangent
vector $V$ to $X$, we can consider the scalar
function $g$ which minimizes the norm
$\norm{V-gX_\xi}$. This function $g$ exists, is
unique and is computed by solving of an elliptic
equation, see Definition \ref{def:FX}. We then
define the seminorm $\tnorm{V}=\norm{V-gX_\xi}$
and consider the distance given by the infimum of
$\int_0^1 \btnorm{\dot X(s)}_{X(s)}\,ds$ over all
paths $X(s)$ joining $X_0$ and $X_1$, that is,
$X(0)=X_0$ and $X(1)=X_1$. The seminorm
$\tnorm{\dott}$ has the property that it vanishes
on the tangent space of all elements that coincide
under relabelings, and, in particular, it implies
that if $X_1$ is a relabeling of $X_0$ then
$d(X_0,X_1)=0$, see Section
\ref{sec:riemann}. With the proper definitions we
find, see Theorem \ref{th:lipstab}, that
$d(\tS_t(X_0),\tS_t(X_1))\leq e^{Ct} d(X_0,X_1)$
for some positive constant $C$, where $\tS_t$
denotes the semigroup that advances the system
\eqref{eq:ode} by a time $t$.  By transfering this
metric to Eulerian variables we finally get a
metric $d_\D$ such that $d_\D(T_t(u,\mu), T_t(\bar
u,\bar\mu))\le e^{Ct}d_\D((u,\mu), (\bar
u,\bar\mu))$, where $T_t$ is the semigroup in
Eulerian variables.

In Section \ref{sec:topology}, we compare the
metric $d_\D$ with other natural topologies. In
particular, in Proposition \ref{prop:convlinf} we
show that if $(u_n, \mu_n)$ converges in the
topology induced by $d_\D$, then $u_n$ converges
in $L^\infty(\Real)$. Furthermore, if $u_n$
converges in $L^\infty(\Real)$ and $u_{x,n}$
converges in $L^2(\Real)$, then the mapping
$u\mapsto (u, u_x^2 dx)$ is continuous on $\D$.

\section{Semi-group of solutions in Lagrangian coordinates}

\subsection{Equivalent system}

The Hunter--Saxton equation equals
\begin{equation}
  \label{eq:hs}
  u_t+uu_x=\frac14(\int_{-\infty}^xu_x^2\,dx-\int_{x}^\infty u_x^2\,dx).
\end{equation}
Formally, the solution satisfies the following
transport equation for the energy density
$u_x^2\,dx$,
\begin{equation}
  \label{eq:transport}
  (u_x^2)_t+(uu_x^2)_x=0
\end{equation}
so that $\int_\Real u_x^2\,dx$ is a preserved
quantity. Next, we rewrite the equation in Lagrangian coordinates. We
introduce the characteristics
\begin{equation*}
  y_t(t,\xi)=u(t,y(t,\xi)).
\end{equation*}
The  Lagrangian velocity  $U$ reads
\begin{equation*}
  U(t,\xi)=u(t,y(t,\xi)).
\end{equation*}
Furthermore, we define the Lagrangian cumulative energy by
\begin{equation*}
  H(t,\xi)=\int_{-\infty}^{y(t,\xi)}u_x^2(t,x)\,dx.
\end{equation*}
From \eqref{eq:hs}, we get that
\begin{equation*}
  U_t=u_t\circ y+y_tu_x\circ y=\frac14\Big(\int_{-\infty}^yu_x^2\,dx-\int_{y}^\infty u_x^2\,dx\Big)=\frac12H(t,\xi)-\frac14H(t,\infty)
\end{equation*}
and
\begin{align*}
  H_t&=\int_{-\infty}^{y(t,\xi)}(u_x^2(t,x))_t\,dx+y_tu_x^2(t,y)\\
  &=\int_{-\infty}^{y(t,\xi)}((u_x^2)_t+(uu_x^2)_x)(t,x)\,dx\\
  &=0
\end{align*}
by \eqref{eq:transport}. Hence, the
Hunter--Saxton equation formally is equivalent to the
following system of ordinary differential
equations:
\begin{subequations}
  \label{eq:sys}
  \begin{align}
    \label{eq:sys1}
    y_t&=U,\\
    \label{eq:sys2}
    U_t&=\frac12H-\frac14H(\infty),\\
    \label{eq:sys3}
    H_t&=0.
  \end{align}
\end{subequations}
We have that $H(\infty)=H_0$ is a constant which
does not depend on time, and global existence of
solutions to \eqref{eq:sys} follows from the
linear nature of the system. There is no exchange
of energy across the characteristics and the
system \eqref{eq:sys} can be solved explicitly,
in contrast with the Camassa--Holm equation where
energy is exchanged across characteristics. We
have
\begin{subequations}
  \label{eq:explsol}
  \begin{align}
    \label{eq:yexpl}
    y(t,\xi)&=(\frac14H(0,\xi)-\frac18H(0,\infty))t^2+U(0,\xi)t+y(0,\xi),\\
    \label{eq:Uexpl}
    U(t,\xi)&=(\frac12H(0,\xi)-\frac14H(0,\infty))t+U(0,\xi),\\
    \label{eq:Uexpl1}
    H(t,\xi)&=H(0,\xi).
  \end{align}
\end{subequations}
Our goal is now to construct a continuous
semigroup of solutions in Eulerian coordinates,
i.e., for the original variable, $u$.  The idea is
to establish a mapping between the variables in
Eulerian and Lagrangian coordinates, and we have
to decide which function space we are going to use
for the solutions of \eqref{eq:sys}. Later, we
will introduce a projection and therefore we need
the framework of Hilbert spaces. A Riemannian
metric also comes from an underlying Hilbert space
structure. Given a natural number $p$, let us
introduce the Banach space (if $p>1$, then
$E^p=H^p(\Real)$)
\begin{equation*}
  E^p=\{f\in L^\infty(\Real) \mid f^{(i)}\in L^2(\Real)\text{ for }i=1,\ldots,p\}
\end{equation*}
and the Hilbert spaces
\begin{equation*}
  H_1^p=H^p(\Real)\times\Real,\quad   H_2^p=H^p(\Real)\times\Real^2.
\end{equation*}
We write $\Real$ as
$\Real=(-\infty,1)\cup(-1,\infty)$ and consider the
corresponding partition of unity $\chi^+$ and
$\chi^-$ (so that $\chi^+$ and $\chi^-\in
C^\infty(\Real)$, $\chi^++\chi^-=1$, $0\leq
\chi^+\leq 1$, $\supp(\chi^+)\subset (-1,\infty)$
and $\supp(\chi^-)\subset (-\infty,1)$). Introduce the linear mapping
$\R_1$ from $H_1^p$ to
$E^p$ defined as
\begin{equation*}
  \xymatrix{(\bar f,a)\ar@{|->}[r]^(0.3){\R_1}&f(\xi)=\bar f(\xi)+a\chi^+(\xi)},
\end{equation*}
and the linear mapping $\R_2$ from $H_2^p$ to $E^p$
defined as
\begin{equation*}
  \xymatrix{(\bar f,a,b)\ar@{|->}[r]^(0.25){\R_2}&f(\xi)=\bar f(\xi)+a\chi^+(\xi)+b\chi^-(\xi).}
\end{equation*}
The mappings $\R_1$ and $\R_2$ are linear,
continuous and injective. Let us introduce $E_1^p$
and $E_2^p$, the images of $H_1^p$ and $H_2^p$ by $\R_1$
and $\R_2$, respectively, that is,
\begin{equation*}
  E_1^p=\R_1(H_1^p)\text{ and }  E_2^p=\R_2(H_2^p).
\end{equation*}
One can check that the mappings $R_1\colon H_1^p\to E_1^p$ and
$R_2\colon H_2^p\to E_2^p$ are homeomorphisms. It follows
that $E_1^p$ can be equipped with two equivalent
norms $\norm{\dott}_E$ and
$\norm{\R_1^{-1}(\dott)}_{H_1^p}$ (and similarly for
$E_2^p$) and, through the mappings $\R_1$ and
$\R_2$, $E_1^p$ and $E_2^p$ can be seen as Hilbert
spaces. We denote
\begin{equation*}
  B^p=E_2^p\times E_2^p\times E_1^p.
\end{equation*}
We will mostly be concerned with the case $p=1$
and to ease the notation, we will not write the
superscript $p$ for $p=1$, that is, $B=B^1$,
$E_j=E_j^1$, etc. In the same way that one proves
that $H^1(\Real)$ is a continous algebra, one
proves the following lemma, which we use later,
\begin{lemma}
  \label{lem:contbialg}
  The space $E$ is a continuous algebra, that is,
  for any $f,g\in E$, then the product $fg$
  belongs to $E$ and there exists constant $C$ such that 
  \begin{equation*}
    \norm{fg}_{E}\leq C\norm{f}_E\norm{g}_E
  \end{equation*}
  for any $f,g\in E$.
\end{lemma}
\begin{definition}
  \label{def:F}
  The set $\F$ consists of the elements
  $(\zeta,U,H)\in B=E_2\times E_2\times E_1$ such
  that
  \\[1mm]
  (i) $(\zeta,U,H)\in(W^{1,\infty})^3$, where $\zeta(\xi)=y(\xi)-\xi$; \\[1mm]
  (ii) $y_\xi\geq0$, $H_\xi\geq0$ and $y_\xi+H_\xi\geq c$, almost everywhere, where  $c$ is a strictly positive constant;\\[1mm]
  (iii) $y_\xi H_\xi=U_\xi^2$ almost everywhere.
\end{definition}
\begin{theorem}
  \label{th:semigroupS}
  The solution of the equivalent system given by
  \eqref{eq:sys} constitutes a semigroup $S_t$ in
  $\F$ which is continuous with respect to the
  $B$-norm. Thus $X(t)=(y(t),U(t),H(t))=S_t(X_0)$
  denotes the solution of \eqref{eq:sys} at time
  $t$ with initial data $X_0$. Moreover, the
  function $\xi\to y(t,\xi)$ is invertible for
  almost every $t$ and we have, for almost every
  $t$, that
  \begin{equation}
    \label{eq:strposyxi}
    y_\xi(t,\xi)>0\text{ for almost every }\xi\in\Real.
  \end{equation}
\end{theorem}
\begin{proof}
  Let
  $(\bar\zeta, \zeta_{\infty},\zeta_{-\infty})$,
  $(\bar U, U_{\infty},U_{-\infty})$ be the
  preimage of $\zeta$ and $U$ by $\R_2$, respectively, and
  $(\bar H, H_{\infty})$ the preimage
  of $H$ by $\R_1$. Inserting these variables into
  \eqref{eq:sys}, we obtain the following linear
  system of equations
  \begin{align*}
    \bar y_t&=\bar U,\\
    \bar U_t&=\frac12\bar H,\\
    \bar H_t&=0,
  \end{align*}
  and
  \begin{align*}
    (y_{\pm\infty})_t&=U_{\pm\infty},\\
    (U_{\pm\infty})_t&=\pm\frac14H_{\infty},\\
    (H_{\pm\infty})_t&=0.
  \end{align*}
  Since it is linear, the system has a global
  solution in $B$, and we have Lipshitz stability with
  respect to the $B$-norm. Again due to the linearity,
  it is clear that the space
  $(W^{1,\infty}(\Real))^3$ is invariant. After
  differentiating \eqref{eq:sys} with respect to
  $\xi$, we obtain
  \begin{subequations}
    \label{eq:dsys}
    \begin{align}
      \label{eq:dsys1}
      y_{\xi t}&=U_\xi,\\
      \label{eq:dsys2}
      U_{\xi t}&=\frac12H_\xi,\\
      \label{eq:dsys3}
      H_{\xi t}&=0.
    \end{align}
  \end{subequations}
  Hence,
  \begin{equation*}
    \frac{d}{dt}(y_\xi H_\xi-U_\xi^2)=0
  \end{equation*}
  so that if the relation
  \begin{equation}
    \label{eq:relvar}
    y_\xi(t,\xi)H_\xi(t,\xi)=U_\xi^2(t,\xi)
  \end{equation}
  holds for $t=0$,  then it holds for all $t$. By assumption, since
  $(y,U,H)_{t=0}\in\F$, we have
  \begin{equation}
    \label{eq:posder}
    (y_\xi+H_\xi)(t,\xi)>0
  \end{equation}
  for $t=0$. By continuity, \eqref{eq:posder} is
  true in a vicinity of $t=0$, and we denote by
  $[0,T)$ the largest interval where it holds. For
  $t\in[0,T)$, it follows from \eqref{eq:relvar}
  that
  \begin{equation}
    \label{eq:posyH}
    y_\xi(t,\xi)\geq0,\quad H_\xi(t,\xi)\geq0,
  \end{equation}
  and
  \begin{equation}
    \label{eq:bdyhbyu}
    \abs{U_\xi}\leq\frac12(y_\xi+H_\xi).
  \end{equation}
  Hence,
  \begin{equation*}
    \frac{d}{dt}(\frac{1}{y_\xi+H_\xi})=\frac{U_\xi}{(y_\xi+H_\xi)^2}
    \leq\frac{1}{2(y_\xi+H_\xi)},
  \end{equation*}
  and, by the Gronwall lemma,
  \begin{equation}
    \label{eq:grinvyH}
    \frac{1}{y_\xi+H_\xi}(t,\xi)\leq\frac{1}{y_\xi+H_\xi}(0,\xi)e^{t/2}
  \end{equation}
  for $t\in[0,T)$. It implies that $T=\infty$ and
  we have proved that $(y(t),U(t),H(t))$ remains
  in $\F$ for all $t$. The proof of statement
  \eqref{eq:strposyxi} goes as in \cite[Lemma
  2.7]{HolRay:07} and we only give here a sketch
  of the argument. Given a fixed $\xi\in\Real$,
  let
  \begin{equation*}
    \N_\xi=\{t\in[0,T]\ |\ y_\xi(t,\xi)=0\}.
  \end{equation*}
  For any $t^*\in\N_\xi$, we have
  \begin{align}
    \label{eq:yxivan1}
    y_\xi(t^*,\xi)&=0,&&\text{from the definition of }t^*,\\
    \label{eq:yxivan2}
    y_{\xi,t}(t^*,\xi)&=0,&&\text{by  \eqref{eq:yxivan1} and \eqref{eq:relvar}, }\\
    \label{eq:yxivan3}
    y_{\xi,tt}(t^*,\xi)&=\frac12
    H_\xi(t^*,\xi)>0,&&\text{by \eqref{eq:yxivan1}
      and \eqref{eq:grinvyH}. }
  \end{align}
  Since the second derivative in time is strictly
  positive, the function $t\to y_\xi(t,\xi)$ is
  strictly positive at least on a small
  neighborhood of $t^*$ excluding $t^*$ where it
  is equal to zero. Note that we can also use the
  explicit formulation given by \eqref{eq:explsol}
  to get the same conclusion. We use Fubini's
  theorem to conclude this argument, see
  \cite{HolRay:07} for the details.
\end{proof}

\subsection{Functional setting in Eulerian
  variables}

Let us define $m=u_{xx}$.
After differentiating
\eqref{eq:hs} twice, we obtain
\begin{equation}
  \label{eq:hsm}
  m_t+um_x+2u_xm=0.
\end{equation}
Note that if we replace $m$ by $u-u_{xx}$, then
\eqref{eq:hsm} will give the Camassa--Holm
equation. For the Camassa--Holm equation there
exists a particular class of solutions that takes the form
\begin{equation*}
  m=\sum_{i=1}^N p_i(t)\delta_{q_i(t)}.
\end{equation*}
Such particular solutions also exist for the
Hunter--Saxton equation, and they correspond to
piecewise linear functions (indeed, $u_{xx}=0$  if $u$ is linear). Let
\begin{equation*}
  y_1(t)=-\frac{t^2}8,\quad U_1(t)=-\frac{t}4,\quad H_1(t)=0,
\end{equation*}
and
\begin{equation*}
  y_2(t)=\frac{t^2}8,\quad U_2(t)=\frac{t}4,\quad H_2(t)=1.
\end{equation*}
Then $(y_1,U_1,H_1)$ and $(y_2,U_2,H_2)$ are
solutions of \eqref{eq:sys} for the total energy
$H(\infty)=1$. One can check that the function $u$
defined as
\begin{equation*}
  u(t,x)=
  \begin{cases}
    U_1(t)&\text{ if }x\leq y_1(t),\\
    \frac{y_1(t)-x}{y_2(t)-y_1(t)} U_1(t)+\frac{x-y_2(t)}{y_2(t)-y_1(t)} U_2(t)&\text{ if }y_1(t)<x\leq y_2(t),\\
    U_2(t)&\text{ if }x>y_2(t),\\
  \end{cases}
\end{equation*}
is a weak solution of \eqref{eq:hs}. At $t=0$, we
have $u(0,x)=0$. However zero is also solution to
\eqref{eq:hs} and therefore, if we want to
construct a semigroup of solution, the function
$u$ at $t=0$ does not provide us with all the
necessary information. We need to know the
location and the amount of energy that has
concentrated on singular set. In the above
example, the whole energy is concentrated at the
origin when $t=0$. The correct space where to
construct global solution of the Hunter--Saxton
equation is given by $\D$ defined as follows.
\begin{definition}
  The set $\D$ consists of all pairs $(u,\mu)$
  such that\\[1mm]
  (i) $u\in E$, $\mu$ is a finite Radon measure;\\[1mm]
  (ii) we have
  \begin{equation}
    \label{eq:condmu}
    \muac=u_x^2dx
  \end{equation}
  where $\muac$ denotes the absolute continuous
  part of $\mu$ with respect to the Lebesgue
  measure.
\end{definition}
We introduce the subset $\F_0$ of $\F$ defined as
follows
\begin{equation}
  \label{eq:defG0}
  \F_0=\{X=(y,U,H)\in\F\mid y+H=\id\}.
\end{equation}
We can define a mapping, denoted $L$, from
$\D$ to $\F_0\subset\F$ as follows.
\begin{definition}
  \label{def:Ldef}
  For any $(u,\mu)$ in $\D$, let
  \begin{subequations}
    \label{eq:Ldef}
    \begin{align}
      \label{eq:Ldef1}
      y(\xi)&=\sup\left\{y\mid \mu((-\infty,y))+y<\xi\right\},\\
      \label{eq:Ldef2}
      H(\xi)&=\xi-y(\xi),\\
      \label{eq:Ldef3}
      U(\xi)&=u\circ{y(\xi)}.
    \end{align}
  \end{subequations}
  Then $X=(\zeta,U,H)\in\F_0$ and we denote by
  $L\colon \D\to\F_0$ the mapping which to any
  $(u,\mu)\in\D$ associates $(\zeta,U,H)\in\F_0$ as
  given by \eqref{eq:Ldef}.
\end{definition}
Thus, from any initial data $(u_0,\mu_0)\in\D$, we
can construct a solution of \eqref{eq:sys} in $\F$
with initial data $X_0=L(u_0,\mu_0)\in\F$. It remains
to go back to the original variables, which is the
purpose of the mapping $M$ defined as follows.
\begin{definition}
  \label{def:Mdef}
  Given any element $X$ in $\F$. Then, the pair
  $(u,\mu)$ defined as follows\footnote{
    The push-forward  of a measure
    $\nu$ by a measurable
    function $f$ is the measure $f_\#\nu$ defined by
    $f_\#\nu(B)=\nu(f\inv(B))$ for all Borel sets $B$.}
  \begin{subequations}
    \label{eq:umudef}
    \begin{align}
      \label{eq:umudef1}
      &u(x)=U(\xi)\text{ for any }\xi\text{ such that  }  x=y(\xi),\\
      \label{eq:umudef2}
      &\mu=y_\#(H_\xi\,d\xi)
    \end{align}
  \end{subequations}
  belongs to $\D$. We denote by
  $M\colon\F\rightarrow\D$ the mapping which to any
  $X$ in $\F$ associates $(u,\mu)$ as given by
  \eqref{eq:umudef}.
\end{definition}
The proofs of the well-posedness of the
Definitions \ref{def:Ldef} and \ref{def:Mdef} are
the same as in \cite[Theorems 3.8 and 3.11]{HolRay:07}.

\subsection{Relabeling symmetry} \label{subsec:label}

When going from Eulerian to Lagrangian coordinates, there
exists an additional degree of freedom which corresponds to
relabeling. Let us explain this schematically. We
consider two elements $X$ and $\bar X$ in $\F$
such that $\bar{X}=X\circ f$, for some function $f$, where $X\circ f$
denotes $(y\circ f,U\circ f,H\circ f)$. The two
element $X$ and $\bar{X}$ correspond to functions
in Eulerian coordinates denoted $u$ and
$\bar{u}$, respectively. We have
\begin{equation*}
  U(\xi)=u\circ y(\xi),\text{ and }\bar{U}(\xi)=\bar{u}\circ \bar{y}(\xi).
\end{equation*}
Then, if $y$ and $\bar{y}$ are invertible, we get
\begin{equation*}
  \bar{u}=\bar{U}\circ \bar{y}^{-1}=U\circ f \circ (y\circ f)^{-1}=U\circ y=u
\end{equation*}
so that $X$ and $\bar{X}$, which may be distinct,
correspond to the same Eulerian configuration. We
can put this statement in a more rigorous
framework by introducing the subgroup $\Gr$ of the
group of homeomorphisms from $\Real$ to $\Real$
defined as
\begin{equation}
  \label{eq:Hcond}
  f-\id\text{ and }f^{-1}-\id\text{ both belong to }W^{1,\infty}(\Real).
\end{equation}
For any $\alpha>1$, we introduce the
subsets $\Gr_\alpha$ of $\Gr$ defined by
\begin{equation*}
  \Gr_\alpha=\{f\in\Gr\mid \norm{f-\id}_{W^{1,\infty}(\Real)}
  +\norm{f^{-1}-\id}_{W^{1,\infty}(\Real)}\leq\alpha\}.
\end{equation*}
The subsets $\Gr_\alpha$ do not possess the group
structure of $G$ but they are closed sets.
\begin{definition}
  \label{def:FGalpha}
  Given $\alpha\geq0$, the set $\F_\alpha$
  (respectively $\G_\alpha$) consists of the
  elements $(\zeta,U,H)\in B=E_2\times E_2\times
  E_1$ such that
  \begin{subequations}
    \label{eq:defFGal}
    \begin{equation}
      \label{eq:defFGal1}
      (\zeta,U,H)\in(W^{1,\infty})^3,
    \end{equation}
    \begin{equation}
      \label{eq:defFGal2}
      y+H\in\Gr_\alpha,
    \end{equation}
    \begin{equation}
      \label{eq:defFGal3}
      y_\xi H_\xi=U_\xi^2 \text{ (respectively $y_\xi H_\xi\geq U_\xi^2$)},
    \end{equation}
    where $\zeta(\xi)=y(\xi)-\xi$.
  \end{subequations}
\end{definition}
We have $\F_\alpha\subset\G_\alpha$. One can check,
using \cite[Lemma 3.2]{HolRay:07}, that
$\F=\bigcup_{\alpha\geq0}\F_\alpha$, and we denote
$\G=\bigcup_{\alpha\geq0}\G_\alpha$. The following
proposition holds.
\begin{proposition}
  (i) The mapping $(X,f)\mapsto \bar{X}$ from
  $\F\times\Gr$ to $\F$ given by $\bar{X}=X\circ f$
  defines an action of the group $\Gr$ on
  $\F$. Hence, we can define the equivalence
  relation on $\F$  by
  \begin{equation*}
    X\sim\bar X\text{ if and only if there exists }f\in\Gr\text{ such that }\bar X=X\circ f,
  \end{equation*}
  and the corresponding quotient is denoted $\quot$.\\[1mm]
  (ii) If $X\sim\bar X$, then $M(X)=M(\bar X)$,
  i.e., the relabeling of an element in $\F$
  corresponds to the same element in $\D$.
\end{proposition}
The proof of this proposition and of the remaining
propositions in this section can be found in
\cite{HolRay:07} with only minor
adaptions. Given $X\in\F$, we denote by $[X]$
the element of $\quot$ which corresponds to the
equivalence class of $X$. We shall see that we can
identify $\quot$ with the subset $\F_0$ of $\F$.
\begin{definition}
  \label{def:Pi}
  We define the projection $\Pi\colon\G\to\G_0$ as
  follows
  \begin{equation*}
    \Pi(X)=X\circ(y+H)^{-1}.
  \end{equation*}
  We have $\Pi(\F)=\F_0$.
\end{definition}
We have the following proposition.
\begin{proposition}
  (i) For $X$ and $\bar X$ in $\F$,
  \begin{equation*}
    \text{$X\sim\bar X$ if and only if  $\Pi(X)=\Pi(\bar X)$.}
  \end{equation*}
  (ii) The injection $X\mapsto [X]$ is a bijection
  from $\F_0$ to $\quot$.
\end{proposition}

\begin{proposition}
  (i) The sets $\D$ and $\F_0$ are in bijection. We
  have
  \begin{equation*}
    M\circ L=\id_{\D}\text{ and }L\circ M|_{\F_0}=\id_{\F_0}.
  \end{equation*}
  (ii) The sets $\D$ and $\quot$ are in bijection.
\end{proposition}

The following proposition says that the solutions
to the system \eqref{eq:sys} are invariant under
relabeling.
\begin{proposition}
  \label{prop:equivariant}
  The mapping $S_t\colon\F\rightarrow\F$ is
  $\Gr$-equivariant, that is,
  \begin{equation}
    \label{eq:Hequivar}
    S_t(X\circ{f})=S_t(X)\circ{f}
  \end{equation}
  for any $X\in\F$ and $f\in\Gr$. This implies that
  \begin{equation*}
    \Pi\circ S_t\circ\Pi=\Pi\circ S_t.
  \end{equation*}
  Hence, we can define a semigroup of solutions on
  $\quot$. It corresponds to the mapping $\tS_t$
  from $\F_0$ to $\F_0$ given by
  \begin{equation}
    \label{eq:deftSt}
    \tS_t=\Pi\circ S_t
  \end{equation}
  which defines a semigroup on $\F_0$.
\end{proposition}

We can rewrite system \eqref{eq:sys} as
\begin{equation}
  \label{eq:sysrew}
  X_t=F(X)
\end{equation}
where $F\colon B\to B$ is given by
\begin{equation}
  \label{eq:defF}
  F(y,U,H)=(U,\frac12 H-\frac14 H(\infty),0).
\end{equation}
Proposition \ref{prop:equivariant} follows from
the fact, which can be verified directly by
looking at \eqref{eq:defF}, that
\begin{equation}
  \label{eq:invfF}
  F(X\circ f)=F(X)\circ f.
\end{equation}
We want to define a distance in $\F_0$ which makes
the semigroup $\tS_t$ Lipschitz continuous.

\section{A Riemannian metric} \label{sec:riemann}

We want to define a mapping $d$ from $\F\times\F$
to $\Real$, which is symmetric and satisfies the
triangle inequality, and such that
\begin{equation}
  \label{eq:sameqre}
  d(X,\bar X)=0\text{ if and only if }X\sim\bar X,
\end{equation}
and
\begin{equation}
  \label{eq:lipdismap}
  d(S_tX,S_ t\bar X)\leq C\, d(X,\bar X),
\end{equation}
because such mapping can in a natural way be used to
define a distance on $\quot$ which also makes the
semigroup of solutions continuous. Since the
stability of the semigroup $S_t$ holds for the
$B$-norm, it is natural to use this norm to
construct the mapping $d$. A natural candidate
would be
\begin{equation*}
  d(X,\bar X)=\inf_{f,\bar f\in\Gr}\norm{X\circ f-\bar X\circ\bar f}_{B},
\end{equation*}
which is likely to fulfill \eqref{eq:sameqre} and
\eqref{eq:lipdismap}. However it does not satisfy
the triangle inequality. Formally, let us explain
our construction, which is inspired by ideas
originating in  Riemannian geometry. Let us think of
$\F$ as a Riemannian manifold embedded in the
Hilbert space $B$. There is a natural scalar
product in the tangent bundle of $T\F$ of $\F$
which is inherited from $B$. We can then define a
distance in $\F$ by considering geodesics, namely,
\begin{equation}
  \label{eq:geodBnorm}
  d(X_0,X_1)=\inf_{X}\int_0^1\norm{\dot X(s)}_{B}\,ds
\end{equation}
for any $X_0,X_1\in\F$ and where the infimum is
taken over all smooth paths $X(s)$ in $\F$ joining
$X_0$ and $X_1$. The distance equals to the
$B$-norm. It makes the semigroup stable but it
clearly separates points which belong to the same
equivalence class and so does not fulfill
\eqref{eq:sameqre}. For a given element $X\in\F$,
we consider the subset $\Gamma\subset\F$ which
corresponds to all relabelings of $X$, that is,
$\Gamma=[X]=\{X\circ f\mid f\in\Gr\}$. If we
substitute in \eqref{eq:geodBnorm} the following
definition
\begin{equation}
  \label{eq:newdist}
  d(X_0,X_1)=\inf_{X}\int_0^1\btnorm{\dot X(s)}_{X(s)}\,ds
\end{equation}
where $\tnorm{\dott}$ is a seminorm in $T\F$ with
the extra property that it vanishes on
$T\Gamma_{X(s)}$, then the property
\eqref{eq:sameqre} will follow in a natural way, and we
expect the stability property \eqref{eq:lipdismap}
to be a consequence of the equivariance of $S_t$,
as stated in Proposition \ref{prop:equivariant}.  We will carry out
the plan next.

Let us first investigate the local structure of $\Gamma$
around $X$. Given a smooth function $g(\xi)$ (one
should actually think of $g$ as an element of
$T\Gr|_{\id}$ ), we consider the curve $f$ in
$\Gr$ given by
\begin{equation*}
  f(\theta,\xi)=\xi+\theta g(\xi).
\end{equation*}
It leads to the curve in $X\circ f(\theta)$ in
$\Gamma$ that we differentiate, and we obtain
\begin{equation*}
  \frac{d}{d\theta}(X\circ f(\theta))=gX_\xi.
\end{equation*}
We now define the subspace $E(X)$ which formally
corresponds to the subspace $T\Gamma_{X}$ of $T\G$.
\begin{definition}
  \label{def:FX}
  Given a fixed element $X\in \G\cap B^2$, we
  consider the subspace $E(X)$ defined as
  \begin{equation*}
    E(X)=\{g(\xi)X_\xi(\xi)
    \mid g\in E_2\},
  \end{equation*}
  where $X_\xi(\xi)=(y_\xi(\xi),U_\xi(\xi),H_\xi(\xi))^T$.
\end{definition}
\begin{lemma}
  \label{lem:accr}
  Given any $X\in B^2$, the bilinear form $a_X$
  defined as 
  \begin{equation*}
    a_X(g,h)=\scal{gX_\xi,hX_\xi}
  \end{equation*}
  is coercive, that is, there exists a constant
  $C>0$ such that
  \begin{equation}
    \label{eq:coerc}
    \frac{1}{C}\norm{g}_{E_2}^2\leq a_X(g,g)=\norm{gX_\xi}_{B}^2
  \end{equation}
  for all $g\in E_2$. Moreover, the constant $C$
  depends only on $\norm{X}_{B^2}$ and
  $\norm{\frac1{y_\xi+H_\xi}}_{L^\infty}$.
\end{lemma}
\begin{proof}
  Given $g\in E_2$, let
  $(\bar g,g_{-\infty},g_{\infty})=\R_2^{-1}(g)$, we
  have the following decomposition, $g=\bar
  g+g_{-\infty}\chi^{-}+g_{\infty}\chi^+$ and, by
  definition,
  \begin{equation*}
    \norm{g}_{E_2}^2=\norm{\bar g}_{H^1}^2+\abs{g_{-\infty}}^2+\abs{g_{\infty}}^2.
  \end{equation*}
  Let us denote $\tilde
  g=g_{-\infty}\chi^{-}+g_{\infty}\chi^+$. Given
  $X\in B^2$, we have
  $\lim_{\xi\to\pm\infty}y_{\xi}(\xi)=1$ and
  $\lim_{\xi\to\pm\infty}(\abs{\zeta_{\xi}}+\abs{U_{\xi}}+\abs{H_\xi})(\xi)=0$. The
  following decomposition hold
  \begin{align*}
    gy_\xi=\bar gy_\xi+\tilde
    g\zeta_\xi+g_{-\infty}\chi^{-}+g_{+\infty}\chi^{+}
  \end{align*}
  so that $\R_2^{-1}(gy_\xi)=(\bar gy_\xi+\tilde
  g\zeta_\xi,g_{-\infty},g_{\infty})$. We have
  also that $\R_2^{-1}(gU_\xi)=(gU_\xi,0,0)$ and
  $\R_1^{-1}(gH_\xi)=(gH_\xi,0)$. Hence,
  \begin{equation*}
    \norm{gX_\xi}_{B}^2=\norm{\bar g y_\xi+\tilde g\zeta_\xi}_{H^1}^2+\norm{gU_\xi}_{H^1}^2+\norm{g H_\xi}_{H^1}^2+\abs{g_{-\infty}}^2+\abs{g_{\infty}}^2.
  \end{equation*}
  Let us prove that 
  \begin{equation}
    \label{eq:coercL2}
    \norm{\bar g}_{L^2}\leq C\norm{gX_\xi}_{B}.
  \end{equation}
  We have
  \begin{multline}
    \label{eq:L2bdbelow}
    \norm{\bar g y_\xi+\tilde
      g\zeta_\xi}_{L^2}^2+\norm{gU_\xi}_{L^2}^2+\norm{g
      H_\xi}_{L^2}^2=\int_{\Real}\bar
    g^2(y_\xi^2+U_\xi^2+H_\xi^2)\,d\xi\\+\int_{\Real}\tilde
    g^2(\zeta_\xi^2+U_\xi^2+H_\xi^2)\,d\xi+2\int_\Real(\bar
    gy_\xi\tilde g\zeta_\xi+\bar gU_\xi\tilde
    gU_\xi+\bar gH_\xi\tilde gH_\xi)\,d\xi
  \end{multline}
  For all $\epsi>0$, we have
  \begin{multline*}
    2\int_\Real(\bar
    gy_\xi\tilde g\zeta_\xi+\bar gU_\xi\tilde
    gU_\xi+\bar gH_\xi\tilde gH_\xi)\,d\xi\geq-\epsi\int_{\Real}\bar
    g^2(y_\xi^2+U_\xi^2+H_\xi^2)\,d\xi\\-\frac1\epsi\int_{\Real}\tilde
    g^2(\zeta_\xi^2+U_\xi^2+H_\xi^2)\,d\xi,
  \end{multline*}
  and, by taking $\epsi$ sufficiently small and
  inserting this inequality into
  \eqref{eq:L2bdbelow}, it yields
  \begin{align*}
    \int_{\Real}\bar
    g^2(y_\xi^2+U_\xi^2+H_\xi^2)\,d\xi&\leq
    C\big(\norm{\bar g y_\xi+\tilde
      g\zeta_\xi}_{L^2}^2+\norm{gU_\xi}_{L^2}^2\\
    &\quad+\norm{g H_\xi}_{L^2}^2+\int_{\Real}\tilde
    g^2(\zeta_\xi^2+U_\xi^2+H_\xi^2)\,d\xi\big)\\
    &\leq C(\norm{gX_\xi}_B^2+\abs{g_{-\infty}}^2+\abs{g_{+\infty}}^2)\\
    &\leq C\norm{gX_\xi}_B^2.
  \end{align*}
  Since $y_\xi^2+U_\xi^2+H_\xi^2(\xi)\geq \frac12
  (y_\xi+H_\xi)^2$, \eqref{eq:coercL2}
  follows. Similarly, by using \eqref{eq:coercL2}
  and a decomposition using $\epsi$ and
  $\frac1\epsi$ as above, one proves that
  \begin{equation*}
    \norm{\bar g_\xi}_{L^2}\leq C\norm{gX_\xi}_B,
  \end{equation*}
  which concludes the proof of the lemma.
\end{proof}

From Lemma \ref{lem:accr} and Lax--Milgram theorem,
we obtain the following definition.
\begin{definition}
  Given any $X\in B^2$ and $V\in B$, there exists
  a unique $g\in E_2$, that we denote $g(X,V)$, such
  that
  \begin{equation}
    \label{eq:defg}
    \scal{gX_\xi,hX_\xi}=\scal{V,hX_\xi}\quad\text{ for all }h\in E_2,
  \end{equation}
  and we have
  \begin{equation*}
    \norm{V-gX_\xi}\leq\norm{X-hX_\xi}\quad\text{ for all }h\in E_2.
  \end{equation*}
\end{definition}
Given $X\in B^2$ and $V\in B$ and $g=g(X,V)$, let
$(\bar g,g_{-\infty},g_\infty)=R_2^{-1}(g)$. When
$V$ is smooth (say $V\in B^2$), one can show that
the following system of equations for $\bar g$,
$g_{-\infty}$ and $g_{\infty}$ is equivalent to
\eqref{eq:defg},
\begin{multline}
  \label{eq:gsys1}
  -\abs{X_\xi}^2\bar g_{\xi\xi}+2(X_{\xi\xi}\cdot
  X_\xi)\bar g_\xi+(\norm{X_\xi}^2+X_\xi\cdot
  X_{\xi\xi\xi})\bar g\\=\bar V-\bar
  V_{\xi\xi}+X_\xi\cdot\big((\id-\partial_\xi^2)\big((g_{-\infty}\chi^-+g_{\infty}\chi^+)[\zeta_\xi,U_\xi,H_\xi]^T\big)\big),
\end{multline}
and
\begin{subequations}
  \label{eq:gsyspm}
  \begin{align}
    \label{eq:gsyspmp}
    (1+\norm{\alpha}_{H^1}^2)g_{\infty}+\scal{\alpha,\beta}_{H^1}g_{-\infty}&=V_\infty-\scal{\bar g X_\xi,\alpha},\\[1mm]
    \label{eq:gsyspmm}
    \scal{\alpha,\beta}_{H^1}g_{\infty}+(1+\norm{\beta}_{H^1}^2)g_{-\infty}&=V_{-\infty}-\scal{\bar
      g X_\xi,\beta},
  \end{align}  
\end{subequations}
where
$\alpha(\xi)=\chi^+(\xi)[\zeta_\xi,U_\xi,H_\xi]^T$
and
$\beta(\xi)=\chi^-(\xi)[\zeta_\xi,U_\xi,H_\xi]^T$
are known functions as they depend only on $X$,
which is given. By Cauchy--Schwarz, the
determinant of system \eqref{eq:gsyspm} for the
unknowns $g_{-\infty}$ and $g_{\infty}$ is strictly
bigger than 1, and therefore we can write
$g_{-\infty}$ and $g_{\infty}$ as functions of
$V$, $X_\xi$ and integral terms which contain
$\bar g$. Since $\abs{X_\xi}^2$ is strictly
bounded away from zero, equation \eqref{eq:gsys1}
for $\bar g$ is elliptic.
\begin{lemma}
  \label{lem:contg}
  The mapping $g:B^2\times B\to E$ is continuous
  and 
  \begin{equation*}
    \norm{g(X_1,V_1)-g(X_0,V_0)}\leq C(\norm{X-\bar X}_{B^2}+\norm{V-\bar V}_B)
  \end{equation*}
  for some constant $C$ which depends only on
  $\norm{V_1}$, $\norm{V_0}$, $\norm{X_1}_{B^2}$,
  $\norm{X_0}_{B^2}$,
  $\norm{(y_{0\xi}+H_{0\xi})^{-1}}_{L^\infty}$,
  $\norm{(y_{1\xi}+H_{1\xi})^{-1}}_{L^\infty}$.
\end{lemma}
\begin{proof}
  From Lemma \ref{lem:contbialg}, it follows that
  $\norm{gX_\xi}_B\leq\norm{g}_{E_2}\norm{X_\xi}_B$
  for any $X\in B^2$ and $g\in E_2$. By
  \eqref{eq:defg} and \eqref{eq:coerc}, we get
  $\norm{g}_{E_2}^2C\leq\norm{V}_B\norm{g}_{E_2}\norm{X_\xi}$
  which implies $\norm{g}_{E_2}\leq C\norm{V}_B$,
  for a constant $C$ which depends only on
  $\norm{X}_{B^2}$. We have, for all $h\in E_2$,
  \begin{multline}
    \label{eq:contproof}
    \scal{(g_1-g_0)X_{1\xi},hX_{1\xi}}=-\scal{g_0(X_{1\xi}-X_{0\xi}),hX_{1\xi}}\\-\scal{g_0X_{0\xi},h(X_{1\xi}-X_{0\xi})}+\scal{V_1-V_0,hX_{1\xi}}+\scal{V_1,h(X_{1\xi}-X_{0\xi})},
  \end{multline}
  which gives
  \begin{equation*}
    \abs{\scal{(g_1-g_0)X_{1\xi},hX_{1\xi}}}\leq C\norm{h}_{E_2}(\norm{V_1-V_0}_{B}+\norm{X_1-X_0}_{B^2}).
  \end{equation*}
  The results follows by taking
  $h=\frac{g_1-g_0}{\norm{g_1-g_0}_{E_2}}$ and
  using \eqref{eq:coerc}.
\end{proof}
We can now define a seminorm on $T\F|_{X}\subset
B$.
\begin{definition}
  Given $X\in B^2$, we define the seminorm
  $\tnorm{\cdot}$ on $B$ as follows: For any element
  $V\in B$, we set
  \begin{equation*}
    \tnorm{V}_{X}=\norm{V-g(X,V)X_\xi}_{B}.
  \end{equation*}
\end{definition}
Using the definition \eqref{eq:newdist} we then
get that
\begin{equation}
  \label{eq:dnulequiv}
  \text{ if }X_0\sim X_1, \text{ then } d(X_0,X_1)=0.
\end{equation}
Indeed, If $X_0\sim X_1$, there exists a function
$f\in\Gr$ such that $X_1=X_0\circ f$. We consider
the path $X(s,\xi)=X_0((1-s)\xi+sf(\xi))$ which joins $X_0$
and $X_1$. We have
\begin{equation*}
  X_s=(f-1)X_{0,\xi}((1-s)\xi+sf(\xi)).
\end{equation*}
Furthermore
\begin{equation*}
  X_\xi=((1-s)\id+sf'(\xi))X_{0,\xi}((1-s)\xi+sf(\xi)).
\end{equation*}
We see that $(1-s)\id+sf'(\xi)\ge \min(\id, f')> 0$.
Thus
\begin{equation*}
  X_s=\frac{(f-1)}{(1-s)\id+sf'(\xi)}X_{\xi},
\end{equation*}
and $X_s\in B$, which implies that $P(X_s)=0$ and therefore
$\tnorm{X_s}_{X(s)}=0$, for all $s\in[0,1]$. Then,
\eqref{eq:dnulequiv} follows from
\eqref{eq:newdist}. In \eqref{eq:newdist}, we
consider the infimum over curves in $\F$. However,
for any $\alpha\geq0$, the set $\F_\alpha$ is not
convex due to the condition \eqref{eq:defFGal3} in
Definition \ref{def:FGalpha}. We relax this
condition and consider instead the set $\G_\alpha$
which is preserved by the semigroup and which is
convex for $\alpha=0$.
\begin{lemma}
  The set $\G_0$ is convex.
\end{lemma}
\begin{proof}
  The set $\Gr_0$ is convex. The condition
  \eqref{eq:defFGal2} implies, for $\alpha=0$,
  that
  \begin{equation}
    \label{eq:yxihxieq1}
    y_\xi+H_\xi=1
  \end{equation}
  which gives $y_\xi^2+2y_\xi
  H_\xi+H_\xi^2=1$. Then the condition
  \eqref{eq:defFGal3} is equivalent to
  \begin{equation}
    \label{eq:convcond}
    y_\xi^2+H_\xi^2+2U_\xi^2\leq 1
  \end{equation}
  which defines a convex set.
\end{proof}
The solution semigroup can be extended to curves
in $\G$. First we define the class of curves we
will be considering.
\begin{definition}
  Given $\alpha\geq0$, we denote by $\C_\alpha$
  the set of curves $X(s)=(\zeta(s),U(s),H(s))$
  where
  \begin{equation*}
    X\colon [0,1]\to \G_\alpha\cap B^2,
  \end{equation*}
  and such that
  \begin{equation*}
    X\in C([0,1],B^2)\quad\text{ and }\quad X_s\in
    C_{\text{{\rm pc}}}([0,1],B)
  \end{equation*}
  where $C_{\text{{\rm pc}}}([0,1],B)$ denotes the set
  of functions from $[0,T]$ to $B$ which are
  piecewise continuous.
\end{definition}
We denote $\C=\bigcup_{\alpha}\C_\alpha$. The
solution operator $S_t$  naturally extends
to curves in $\C$.
\begin{lemma}
  \label{lem:curvesol}
  For any initial curve $X_0\in\C$, there exists a
  solution curve
  $X\colon [0,1]\times \Real_+\to B^2$ such that\\
  (i) $X(s,0)=X_0(s)$;\\
  (ii) for each fixed $t\in\Real_+$,
  $X(\dott,t)\colon[0,1]\to B^2$ belongs to
  $\C$; \\
  (iii) for each fixed $s\in[0,1]$,
  $X(s,\dott)\colon\Real_+\to B^2$ is a solution of
  \eqref{eq:sys} with initial data $X_0(s)$.\\
  Moreover, we have
  \begin{equation}
    \label{eq:lowbdal}
    (y+H)(t,\dott)\in\Gr_{\alpha(t)}\text{ with }\alpha(t)\leq e^{Ct}
  \end{equation}
  for some constant $C$.
\end{lemma}
\begin{proof}
  The proof follows as the proof of
  Theorem \ref{th:semigroupS}. We use a fixed
  point argument, for $T$ small enough, on the set
  $C([0,T],\bar\C)$ where $\bar\C$ is the Banach
  space of curves with piecewise constant
  derivatives, i.e.,
  \begin{equation*}
    \bar\C=\{X\in C([0,1],B^2)\mid X_s\in C([s_i,s_{i+1}],B),\ i=1,\ldots,n\}
  \end{equation*}
  where the sequence $0=s_1\leq\cdots\leq s_n=1$
  is chosen such that $X_{0}\in\bar C$. We then
  extend the solution globally in time and obtain
  \eqref{eq:lowbdal} as in the proof of Theorem
  \ref{th:semigroupS}.
\end{proof}
Define a metric on $\G_0$ as follows.
\begin{definition}
  For two elements $X_0,X_1\in \G_0\cap B^2$, we
  define
  \begin{equation}
    \label{eq:defdist}
    d(X_0,X_1)=\inf_{X\in\C_0}\int_{0}^1\tnorm{X_s(s)}_{X(s)}\,ds.
  \end{equation}
\end{definition}
Note that the definition is well-posed because
$\C_0$ is nonempty since, as $\G_0\cap B^2$ is
convex, we can always join two elements in $C_0$
by a straight line.
\begin{lemma}
  \label{lem:proofdist}
  The mapping $d\colon \G_0\times\G_0\to\Real_+$ is a
  distance on $\G_0\cap B^2$.
\end{lemma}
\begin{proof}
  Let us first prove that $d(X_0,X_1)=0$ implies
  $X_0=X_1$. For any $\epsi\geq0$, we consider
  $X\in\C_0$ such that
  \begin{equation}
    \label{eq:deqepsi}
    \int_0^1\tnorm{X_s}_{X(s)}\,ds\leq\epsi.
  \end{equation}
  Since $y(s,\xi)+H(s,\xi)=\xi$ for all $\xi$, we
  get
  \begin{equation*}
    y_s+H_s=0,\quad\text{ and }\quad y_\xi+H_\xi=1.
  \end{equation*}
  We consider the orthogonal decomposition of
  $X_s$, i.e.,
  \begin{equation}
    \label{eq:decXsg}
    X_s(s,\xi)=g(s,\xi)X_\xi(s,\xi)+R(s,\xi).
  \end{equation}
  It follows, by adding the first and third
  components in \eqref{eq:decXsg}, that
  \begin{equation*}
    0=y_s+H_s=g(s,\xi)(y_\xi+H_\xi)+R_1+R_3=g(s,\xi)+R_1+R_3
  \end{equation*}
  (where $R_1$ and $R_3$ denotes the first and
  third components of $R$) and therefore
  \begin{equation}
    \label{eq:decR1R2}
    g(s,\xi)=R_1(s,\xi)+R_3(s,\xi).
  \end{equation}
  Since, in a Euclidean space the shortest path
  between two points is a straight line, we have
  \begin{equation}
    \label{eq:eucli}
    \norm{X_1-X_0}_{L^\infty(\Real)}\leq\int_0^1\norm{X_s(s,\dott)}_{L^\infty}\,ds.
  \end{equation}
  From the definition of $\G_0$, it follows that
  $y_\xi$, $H_\xi$ and $U_\xi$ are bounded by one
  in $L^\infty(\Real)$, see \eqref{eq:convcond}.  Therefore,
  \eqref{eq:eucli} and
  \eqref{eq:decXsg} imply
  \begin{align*}
    \norm{X_1-X_0}_{L^\infty(\Real)}
    &\leq\int_0^1(\norm{g(s,\dott)}_{L^\infty}+\norm{R(s,\dott)}_{L^\infty})\,ds\\
    &\leq
    2\int_0^1\norm{R(s,\dott)}_{L^\infty}\,ds\quad\text{ (by \eqref{eq:decR1R2})}\\
    &\leq
    2\int_0^1\norm{R(s,\dott)}_B\,ds=2\int_0^1\tnorm{X(s,\dott)}_{X(s)}\,ds
    \leq\epsi.
  \end{align*}
  Since $\epsi$ is arbitrary, it follows that
  $X_1=X_0$. The triangle inequality is obtained
  by patching two curves together and
  reparametrizing them while the symmetry of $d$
  is also obtained by reparametrization. Both
  proofs are somehow standard.
\end{proof}
On $\G_0$, the distance $d$ is weaker than the
$B$-norm as the next lemma shows.
\begin{lemma}
  \label{lem:bddB}
  For any $X_0,X_1\in\G_0\cap B^2$, we have
  \begin{equation}
    \label{eq:ineqdi}
    d(X_0,X_1)\leq\norm{X_1-X_0}.
  \end{equation}
\end{lemma}
\begin{proof}
  Consider $X\in\C_0$ defined as
  follows
  \begin{equation*}
    X(s)=(1-s)X_0+sX_1.
  \end{equation*}
  We have
  \begin{equation*}
    d(X_0,X_0)\leq \int_0^1\tnorm{X_s(s)}_{X(s)}\,ds\leq\int_0^1\norm{X_s(s)}_{X(s)}\,ds=\norm{X_1-X_0}
  \end{equation*}
  because
  $\tnorm{X_s}=\norm{P(X_s)}\leq\norm{X_s}$ as $P$
  is an orthogonal projection.
\end{proof}

\begin{definition}
  \label{def:dist2}
  For two elements $X_0,X_1\in \G_0$,
  we define
  \begin{equation}
    \label{eq:defdist2}
    d(X_0,X_1)=\lim_{n\to\infty}d(X_0^n,X_1^n)
  \end{equation}
  for any sequences $X_0^n$ and $X_1^n$ in $\G_0\cap B^2$
  which converge in $B$ to $X_0$ and $X_1$,
  respectively.
\end{definition}
Definition \ref{def:dist2} is well-posed thanks to
\eqref{eq:ineqdi}. The mapping $S_t$ maps $\F$ to
$\F$, and we can formally define what is called in
differential geometry the tangent map of $S_t$,
$TS_t$, which is a mapping $T\F_{X}$ to
$T\F_{S_tX}$. The following theorem expresses the
fact that $TS_t$ is uniformly continuous (in time)
with respect to the seminorm $\tnorm{\cdot}$.
\begin{theorem}
  \label{th:lipbdnorm}
  There exists a constant $C$ such that, for any
  initial curve $X_0(s,\xi)\in\C_0$, if we
  consider the curve solution $X(t,s,\xi)$ with
  initial data $X_0(s,\xi)$ given by Lemma
  \ref{lem:curvesol}, we have
  \begin{equation}
    \label{eq:lipsemnorm}
    \tnorm{X_s(s,t)}_{X(s,t)}\leq e^{Ct}\tnorm{X_s(s,0)}_{X(s,0)}.
  \end{equation}
\end{theorem}
\begin{proof}
  We rewrite the system
  \begin{equation}
    \label{eq:simpsys}
    X_t=F(X)
  \end{equation}
  where $F$ is given by \eqref{eq:defF}. The
  mapping $F$ is linear and therefore
  differentiable and we have, for any $X,\bar X\in
  B$,
  \begin{equation}
    \label{eq:defDF}
    DF[X](\bar X)=F(\bar X)
  \end{equation}
  where $DF[X]$ denotes the diffential of $F$ at
  $X$. For $X\in B^2$, since $X_\xi\in H^1$, we
  have $\lim_{\xi\to\infty} H_\xi(\xi)=0$ and one
  can then check directly that, for any $g\in
  E_2$,
  \begin{equation}
    \label{eq:fundrel}
    DF[X](g(\xi)X_\xi(\xi))=g(\xi)(DF[X](X_\xi(\xi))).
  \end{equation}
  However, the simplicity of system \eqref{eq:sys}
  may hide the more fundamental nature of relation
  \eqref{eq:fundrel}, which in fact corresponds to
  the infinitesimal version of the equivariance
  property of $F$ stated in
  \eqref{eq:invfF}. Indeed, given a smooth
  function $g$, we consider the family of
  diffeomorphisms parametrized by $\theta$ given by
  $f^\theta(\xi)=\xi+\theta g(\xi)$. The
  equivariance property \eqref{eq:invfF} of $F$
  gives
  \begin{equation*}
    F(X\circ f^\theta)=F(X)\circ f^\theta,
  \end{equation*}
  which after differentiation by $\theta$ and
  taking the value at $\theta=0$ yields
  \eqref{eq:fundrel}. After differentiating
  \eqref{eq:simpsys} with respect to $s$, we get
  \begin{equation}
    \label{eq:difsympsyswrs}
    X_{st}=DF[X](X_s)
  \end{equation}
  while differentiating it with respect to $\xi$
  yields
  \begin{equation}
    \label{eq:difsympsyswrxi}
    X_{\xi t}=DF[X](X_\xi).
  \end{equation}
  We consider the decomposition of $X_s$ given by
  \begin{equation}
    \label{eq:decXS}
    X_s=g(X,V)X_\xi+R.
  \end{equation}
  Since, for every $s\in[0,1]$, $X_s\in
  C^1([0,T],B)$, $X_\xi\in C^1([0,T],B^2)$ and
  \eqref{eq:lowbdal} holds, we can use Lemma
  \ref{lem:contg} to prove that $g\in
  C^1([0,T],E_2)$, for any $s\in[0,1]$. By
  differentiating
  \begin{equation*}
    \scal{gX_\xi,hX_\xi}=\scal{X_s,hX_\xi}
  \end{equation*}
  we obtain that $g_t$ is defined as the unique
  element in $E_2$ such that
  \begin{multline*}
    \scal{g_t X_\xi,hX_\xi}=\scal{X_{st},hX_\xi}+\scal{X_s,hX_{\xi t}}-\scal{gX_{\xi t},hX_\xi}-\scal{gX_\xi,hX_{\xi t}}
  \end{multline*}
  for all $h\in E_2$.  We differentiate
  \eqref{eq:decXS} and get
  \begin{equation*}
    X_{st}=g_t X_\xi+g X_{\xi t}+R_t.
  \end{equation*}
  After using \eqref{eq:difsympsyswrs} and
  \eqref{eq:difsympsyswrxi}, it yields
  \begin{equation*}
    DF[X](X_s)=g_tX_\xi+g(DF[X](X_\xi))+R_t.
  \end{equation*}
  Using \eqref{eq:fundrel}, this identity rewrites
  \begin{equation*}
    R_t=DF[X](X_s-gX_\xi)-g_tX_\xi
  \end{equation*}
  or
  \begin{equation}
    \label{eq:exprRt}
    R_t=DF[X]R-g_tX_\xi.
  \end{equation}
  We take the scalar product of $R_t$ and, since
  $g_tX_\xi$ and $R$ are orthogonal, we obtain
  \begin{align}
    \notag \scal{R_t,R}&=\scal{DF[X](R),R}\\
    \notag
    &\leq\norm{DF[X](R)}\norm{R}\\
    \label{eq:bdscaprRRt}
    &\leq C\norm{R}^2
  \end{align}
  because the mapping $DF[X]\colon B\to B$ is uniformly
  bounded, see \eqref{eq:defDF}. Thus,
  \eqref{eq:bdscaprRRt} yieds
  \begin{equation*}
    \frac{d}{dt}\norm{R}^2\leq C\norm{R}^2.
  \end{equation*}
  By Gronwall's inequality, it implies
  \begin{equation*}
    \tnorm{X_s(t)}=\norm{R(t)}\leq\norm{R(0)} e^{Ct}=\tnorm{X_s(0)} e^{Ct}.
  \end{equation*}
\end{proof}

\begin{theorem}
  \label{th:lipstab}
  The semigroup $\tS_t\colon\G_0\to\G_0$ is Lipschitz
  continuous with respect to the metric $d$. We
  have, for some constant $C$,
  \begin{equation}
    \label{eq:lipstab}
    d(\tS_t(X_0),\tS_t(X_1))\leq e^{Ct} d(X_0,X_1)
  \end{equation}
  for all $X_0,X_1\in\G_0$.
\end{theorem}
\begin{proof}
  We consider first initial conditions
  $X_0,X_1\in\F_0$. There exists a curve $X(s)$ in
  $\C_0$ such that
  \begin{equation*}
    \int_{0}^1\tnorm{X_s(s)}_{X(s)}\,ds\leq d(X_0,X_1)+\epsi.
  \end{equation*}
  We consider the corresponding solution given by
  Lemma \ref{lem:curvesol}, that we simply denote
  $X(s,t)$. By Theorem \ref{th:lipbdnorm}, we have
  \begin{equation}
    \label{eq:XsX0bd}
    \tnorm{X_s(s,t)}_{X(s,t)}\leq e^{Ct}\tnorm{X_s(s,0)}_{X(s,0)}.
  \end{equation}
  Given a time $T$, we consider the projection of
  the curve $X(s,T,\dott)$ on $\G_0$, that we denote
  $\bar X(s,\xi)$, which is given by
  \begin{equation*}
    \bar X(s,\dott)=\Pi(X(s,T,\dott)).
  \end{equation*}
  We denote by $f(s,t,\xi)$, the inverse of
  $(y+H)(s,t,\xi)$ with respect to $\xi$, which is
  allways well-defined and bounded as
  $(y+H)(s,t,\dott)\in\Gr_\alpha$ for some
  $\alpha\leq e^{Ct}$, see Lemma
  \ref{lem:curvesol}. The definition of $\Pi$
  gives
  \begin{equation*}
    \bar X(s,\xi)=X(s,T,f(s,T,\xi)).
  \end{equation*}
  We have $X(0,\dott)=\tS_TX_0$ and
  $X(1,\dott)=\tS_TX_1$ and the curve $\bar X$
  belongs to $\C_0$. We have
  \begin{equation}
    \label{eq:barXs}
    \bar X_s(s,\xi)=X_s(s,T,f)+f_s X_\xi(s,T,f)
  \end{equation}
  and
  \begin{equation}
    \label{eq:barXxi}
    \bar X_\xi(s,\xi)=f_\xi X_s(s,T,f).
  \end{equation}
  We consider decomposition of $X_s$ given by
  \begin{equation}
    \label{eq:decXS2}
    X_s(s,T,\xi)=g(s,T,\xi)X_\xi(s,T,\xi)+R(s,T,\xi).
  \end{equation}
  where
  $g(s,T,\cdot)=g(X(s,T,\cdot),X_s(s,T,\cdot))$.
  Combining \eqref{eq:barXs}, \eqref{eq:barXxi}
  and \eqref{eq:decXS2}, we end up with
  \begin{equation}
    \label{eq:predecbarXs}
    \bar X_s(s,\xi)=\left(\frac{g(s,T,f(\xi))}{f_\xi(s,T,\xi)}+f_s(s,T,\xi)\right)\bar X_\xi(s,\xi)+R(s,t,f(s,T,\xi)).
  \end{equation}
  Hence,
  \begin{equation}
    \label{eq:bdtnoXs}
    \tnorm{\bar X_s(s,\xi)}\leq\norm{R(s,t,f(s,T,\xi))}.
  \end{equation}
  Let us prove that
  \begin{equation}
    \label{eq:bdRcompf}
    \norm{R(s,T,f(s,T,\xi))}\leq
    e^{Ct}\norm{R(s,T,\xi)}
  \end{equation}
  for some constant $C$. We have to prove
  that for any $g\in E_2$, we have
  \begin{equation}
    \label{eq:compnorel}
    \norm{g\circ f}_{E_2}\leq
    e^{Ct}\norm{g}_{E_2}.
  \end{equation}
  We have
  \begin{equation}
    \label{eq:bdlinftyg}
    \norm{g\circ
      f}_{L^\infty(\Real)}\leq\norm{g}_{L^\infty(\Real)}
  \end{equation}
  and
  \begin{align}
    \notag \norm{(g\circ
      f)_\xi}_{L^2(\Real)}^2&=\int_\Real
    (g_\xi\circ f)^2 f_\xi^2\,d\xi\\
    \notag
    &\leq \norm{f_\xi}_{L^\infty(\Real)}\int_\Real
    (g_\xi\circ f)^2 f_\xi\,d\xi\\
    \label{eq:bdgcircf}
    &=\norm{f_\xi}_{L^\infty(\Real)}\norm{g}_{L^2(\Real)}^2.
  \end{align}
  Since $f_\xi=\frac1{y_\xi+H_\xi}\circ f$, we
  have $\norm{f_\xi}_{L^\infty(\Real)}\leq e^{CT}$
  by \eqref{eq:lowbdal}, as
  $(y_\xi+H_\xi)(s,0,\xi)=1$ for all $\xi$, and
  \eqref{eq:bdlinftyg}, \eqref{eq:bdgcircf} imply
  \eqref{eq:compnorel}. Using \eqref{eq:bdRcompf},
  it follows from \eqref{eq:bdtnoXs} that
  \begin{equation}
    \label{eq:barXsbd}
    \tnorm{\bar X_s(s,\xi)}\leq C\tnorm{X_s(s,T,\xi)}
  \end{equation}
  because
  $\norm{R(s,T,\xi)}=\tnorm{X_s(s,T,\xi)}$. Hence,
  we finally get
  \begin{align*}
    d(\tS_tX_0,\tS_tX_1)&\leq\int_{0}^1\tnorm{\bar
      X_s}_{\bar X(s)}\,ds\\
    &\leq e^{CT}\int_{0}^1\tnorm{ X_s(s,T)}_{X(s,T)}\,ds\ \text{ (by \eqref{eq:barXsbd})}\\
    &\leq e^{2CT}\int_{0}^1\tnorm{ X_s(s,0)}_{X(s,0)}\,ds\ \text{ (by \eqref{eq:XsX0bd})}\\
    &\leq e^{2CT}(d(X_0,X_1)+\epsi)
  \end{align*}
  which implies \eqref{eq:lipstab} as $\epsi$ is
  arbitrary. To extend this result to any
  $X_0,X_1\in\G_0$, we use the fact that the
  mapping $\tS_t$ is continuous with respect to
  the $B$-norm (Lemma \ref{lem:tscont}) and
  $\G_0\cap B^2$ is dense in $\G_0$ (Lemma
  \ref{lem:g0capbdense}).
\end{proof}
\begin{lemma}
  \label{lem:tscont}
  The mapping $\Pi\colon\F_\alpha\to\F_0$ is continuous
  with respect to the $B$-norm. It follows that
  $\tS_t$ is a continuous semigroup with respect
  to the $B$-norm.
\end{lemma}
\begin{proof}
  The proof of the continuity of $\Pi$ is the same
  as in in \cite[Lemma 3.5]{HolRay:07}. The
  continuity of $\tS_t$ then follows from
  \eqref{eq:deftSt} and the fact that
  $S_t\colon\F_0\to\F_{\alpha(t)}$ for $\alpha\leq
  e^{Ct}$.
\end{proof}
\begin{lemma}
  \label{lem:g0capbdense}
  The set $\G_0\cap B^2$ is dense in $\G_0$.
\end{lemma}
\begin{proof}
  Given $X_0\in\G_0$, we first assume that
  $X_{0,\xi}$ has compact support. We consider a
  mollifier $\rho^\epsi$. Given $X\in\G_0$, we
  consider the approximation
  $X^\epsi=X\star\rho^\epsi=(\zeta\star\rho^\epsi,U\star\rho^\epsi,H\star\rho^\epsi)$. By the
  Jensen inequality, since $\rho^\epsi\geq0$ and
  $\int_{\Real}\rho^\epsi(\eta)\,d\eta=1$, we have
  \begin{equation*}
    \left(\int_{\Real} \zeta_\xi(\xi-\eta)\rho^\epsi(\eta)\,d\eta\right)^2\leq\int_{\Real} \zeta_\xi(\xi-\eta)^2\rho^\epsi(\eta)\,d\eta
  \end{equation*}
  and similar inequalities for $U_\xi$ and
  $H_\xi$. Hence, since $X$ satisfies
  \eqref{eq:convcond},
  \begin{align*}
    ((y_\xi^\epsi)^2+(H_\xi^\epsi)^2+2(U_\xi^\epsi)^2)(\xi)&\leq\int_\Real
    ((y_\xi)^2+(H_\xi)^2+2(U_\xi)^2)(\xi-\eta)\rho^\epsi(\eta)\,d\eta\\
    &\leq\int_\Real\rho^\epsi(\eta)\,d\eta=1,
  \end{align*}
  and $X^\epsi$ also satisfies
  \eqref{eq:convcond}. Since $y+H=\id$, we have
  \begin{equation*}
    y^\epsi+H^\epsi=\int_{\Real} (\xi-\eta)\rho^\epsi(\eta)\,d\eta=\xi
  \end{equation*}
  (we consider an even mollifier) and $X^\epsi$
  satisfies \eqref{eq:defFGal2} for
  $\alpha=0$. Since $X_\xi$ has a compact support,
  which we denote $K$, $X(\xi)$ is constant for
  $\xi\in K^c$ and $X^\epsi=X$ on a the complement
  of a compact neighborhood of $K$, for $\epsi$
  small enough. Since $X^\epsi\to X$ on any
  compact set, it follows that $X^\epsi\to X$ in
  $L^\infty(\Real)$. By the standard convergence
  properties of approximating sequences, we have
  $X_\xi^\epsi\to X_\xi$ in $L^2(\Real)$ so that,
  finally, $X^\epsi\to X$ in $B$. Let us now
  consider the case where $X\in\G_0$ does not have
  a compact support. For any integer $n$, we
  define $X^n\in\G_0$ as follows
  \begin{equation*}
    X^n(\xi)=
    \begin{cases}
      X^n(-n)&\text{ if } \xi\leq-n,\\
      X^n(\xi)&\text{ if }-n<\xi<n,\\
      X^n(n)&\text{ if }\xi\geq n.
    \end{cases}
  \end{equation*}
  We have
  \begin{equation*}
    X^n_\xi=
    \begin{cases}
      X_\xi&\text{ if }\xi\in(-n,n),\\
      0&\text{ otherwise, }
    \end{cases}
  \end{equation*}
  so that $X_\epsi^n$ has a compact support and
  the condition \eqref{eq:defFGal3} is
  satisfied. Since $X\in B$, we have
  $\lim_{\eta\to\pm\infty}X(\xi)=X(\pm\infty)$ and
  $X^n$ tends to $X$ in $L^\infty(\Real)$. Since
  $X^n_\xi$ is a cut-off of $X_\xi$ with a growing
  support, $X_\xi^n$ tends to $X_\xi$ in
  $L^2(\Real)$. Therefore $X^n$ tends to $X$ in
  $B$ and we have proved that the functions
  $X\in\G_0$ such that $X_\xi$ has compact support
  are dense in $\G_0$.
\end{proof}

\section{Semi-group of solutions in Eulerian coordinates}\label{sec:euler}

We now return to the Eulerian variables.
\begin{definition}\label{def:semiE}
  Let
  \begin{equation}
    \label{eq:semiE}
    T_t=MS_tL\colon \D\to \D.
  \end{equation}
\end{definition}
Next we show that $T_t$ is a Lipschitz continuous semigroup
by introducing a metric on $\D$.

Using the bijection $L$ we
can transport the topology from $\F_0$ to $\D$.
\begin{definition}\label{def:metricE}
  Define the metric $d_\D\colon \D\times\D\to [0,\infty)$ by
  \begin{equation}
    \label{eq:metriE}
    d_\D((u,\mu),(\bar u,\bar\mu))=d(L(u,\mu),L(\bar u,\bar\mu)).
  \end{equation}
\end{definition}
The final result in Eulerian variables reads as follows.
\begin{theorem}
  We have that $(T_t,d_\D)$ is a continuous semigroup on $\D$.
\end{theorem}
\begin{proof}
  We have the following calculation
  \begin{align*}
    d_\D(T_t(u,\mu), T_t(\bar u,\bar\mu))&
    = d(L(T_t(u,\mu)), L(T_t(\bar u,\bar\mu))) \\
    &=d(LT_tML(u,\mu), LT_tML(\bar u,\bar\mu)) \\
    &=d(S_tL(u,\mu), S_tL(\bar u,\bar\mu)) \\
    &\le e^{Ct}d(L(u,\mu), L(\bar u,\bar\mu)) \\
    &=e^{Ct}d_\D((u,\mu), (\bar u,\bar\mu)).
  \end{align*}
\end{proof}

By a weak solution of \eqref{eq:hs0} we mean the
following.
\begin{definition}
  Let $u\colon \Real\times \Real \to \Real$ that satisfies: \\
  (i) $u\in C([0,\infty);L^\infty(\Real))$ and $u_x\in L^\infty([0,\infty);L^2(\Real))$; \\
  (ii) the equation
  \begin{equation}
    \label{eq:weakform}
    \iint_{[0,\infty)\times\Real}\big(u \phi_t-(uu_x-V)\phi \big)\, dxdt
    = \int_\Real u_0\phi|_{t=0}\, dx
  \end{equation}
  holds for all $\phi\in C^\infty_0([0,\infty)\times\Real)$. Here
  $V(t,x)=\frac14(\int_{-\infty}^xu_x^2\,dx-\int_{x}^\infty
  u_x^2\,dx)$ is in $L^\infty([0,\infty);L^\infty(\Real))$. Then we say that $u$ is a weak global conservative solution of the Hunter--Saxton equation \eqref{eq:hs0}.
\end{definition}

\begin{theorem} Given any initial condition
  $(u_0,\mu_0)\in\D$, we denote
  $(u,\mu)(t)=T_t(u_0,\mu_0)$. Then, $u(t,x)$ is a
  global solution of the Hunter--Saxton equation.
\end{theorem}

\begin{proof}
  After making the change of variables
  $x=y(t,\xi)$, we get, on the one  hand,
  \begin{align}
    \notag \iint_{[0,\infty)\times\Real} u \phi_t\,
    dxdt&=\iint_{[0,\infty)\times\Real} u(t,y(t,\xi)) \phi_t(t,y(t,\xi))y_\xi(t,\xi)\, d\xi dt\\
    \notag
    &=\iint_{[0,\infty)\times\Real}U(\phi(t,y)_t-y_t\phi_x(t,y))y_\xi\, d\xi dt\\
    \notag
    &=-\iint_{[0,\infty)\times\Real} (U_ty_\xi+y_{\xi,t} U)\phi(t,y)\, d\xi dt\\
    \notag &\quad -\iint_{[0,\infty)\times\Real}
    U^2\phi(t,y)_\xi\, d\xi dt\\ \notag
    &\quad+\int_\Real
    (Uy_\xi)(0,\xi)\phi(0,y(0,\xi))y_\xi(0,\xi)\,d\xi\\
    \notag
    &=-\iint_{[0,\infty)\times\Real} \big((\frac12H-\frac14H(\infty))y_\xi\big)\phi(t,y)\, d\xi dt\\
    \notag
    &\quad+\iint_{[0,\infty)\times\Real} (U_\xi U)\phi(t,y)\, d\xi dt\\
    \label{eq:solpf1}
    &\quad+\int_\Real
    u(0,x)\phi(0,x)\,dx,
  \end{align}
  and, on the other hand,
  \begin{align}
    \notag \iint_{[0,\infty)\times\Real}\big(uu_x&-V
    \big)\phi\, dxdt\\
    \notag &=\iint_{[0,\infty)\times\Real}\big(uu_x-V
    \big)(t,y)\phi(t,y)y_\xi\,
    d\xi dt\\
    \notag &=\iint_{[0,\infty)\times\Real}\big(UU_\xi
    \big)\phi(t,y)\,
    d\xi dt\\
    \label{eq:solpf2}
    &\quad-\iint_{[0,\infty)\times\Real}V(t,y)\phi(t,y)y_\xi\,
    d\xi dt.
  \end{align}
  By using \eqref{eq:defFGal3} and the fact that
  $U_\xi=u_x\circ y y_\xi$, we get
  \begin{equation}
    \label{eq:h1bd0}
    \int_{\Real}u_x^2\,dx=\int_{\Real}u_x^2\circ  y y_\xi\,d\xi
    =\int_{\{\xi\in\Real\mid y_\xi(t,\xi)>0\}}\frac{U_\xi^2}{y_\xi}\,d\xi=\int_{\{\xi\in\Real\mid y_\xi(t,\xi)>0\}}H_\xi\,d\xi.
  \end{equation}
  The statement \eqref{eq:strposyxi} implies that,
  for almost every $t\in\Real$, the set $\{\xi\in\Real\mid
  y_\xi(t,\xi)>0\}$ is of full measure and
  therefore \eqref{eq:h1bd0} yields
  \begin{equation}
    \label{eq:h1bd}
    \int_{\Real}u_x^2\,dx=\int_{\Real}H_\xi\,d\xi=H(\infty),
  \end{equation}
  for almost every $t\in\Real$. Similarly, for almost
  every $t\in\Real$, we get
  \begin{align}
    \notag
    V(t,y(t,\xi))&=\frac12\int_{\infty}^{y(t,\xi)}u_x^2\,dx-\frac14\int_{\Real}u_x^2\,dx\\
    \notag
    &=\frac12\int_{\infty}^{\xi}u_x^2(t,y(t,\xi))y_\xi(t,\xi)\,dx-\frac14 H(\infty)\\
    \notag
    &=\frac12\int_{\infty}^{\xi}H_\xi(t,\xi)\,d\xi-\frac14H(\infty)\\
    \label{eq:solpf3}
    &=\frac12H(t,\xi)-\frac14 H(\infty).
  \end{align}
  After gathering \eqref{eq:solpf1},
  \eqref{eq:solpf2} and \eqref{eq:solpf3}, we
  obtain that $u$ is a weak solution of the
  Hunter--Saxton equation. It follows from
  \eqref{eq:h1bd0} that
  \begin{equation*}
    \int_{\Real}u_x^2(t,x)\,dx\leq H(t,\infty)=H(0,\infty)=\mu_{0}(\Real)
  \end{equation*}
  so that $u_x\in L^\infty(\Real,L^2(\Real))$. By
  construction of the semigroup $T_t$, we know
  that $(u,\mu)(t)\in C(\Real,\D)$ where $\D$ is
  equipped by the metric $d_\D$. Proposition
  \ref{prop:convlinf} below then implies that
  $u\in C(\Real,L^\infty(\Real))$.
\end{proof}

\section{The topology induced by the metric $d_D$} \label{sec:topology}

\begin{proposition}
  \label{prop:convE}
  The mapping
  \begin{equation*}
    u\mapsto(u,u_x^2\,dx)
  \end{equation*}
  is continuous from $E_2$ into $\D$. In other
  words, given a sequence $u_n\in E_2$ converging to
  $u$ in $E_2$, that is,
  \begin{equation*}
    u_n\to u\text{ in }L^\infty(\Real)\text{ and }u_{n,x}\to u_x\text{ in }L^2(\Real),
  \end{equation*}
  then $(u_n,u_{nx}^2\,dx)$ converges to $(u,u_x^2\,dx)$
  in $\D$.
\end{proposition}
\begin{proof}
  Let $X_n=(y_n,U_n,H_n)=L(u_n,u_{nx}^2\,dx)$ and
  $X=(y,U,H)=L(u,u_x^2\,dx)$, see
  \eqref{eq:Ldef}. Following the proof of
  \cite[Proposition 5.1]{HolRay:07}, one can prove
  that
  \begin{equation*}
    X_n\to X\text{ in }B.
  \end{equation*}
  Hence, by \eqref{eq:ineqdi} in Lemma
  \ref{lem:bddB}, we get that
  $\lim_{n\to\infty}d(X_n,X)=0$ and therefore
  \begin{equation*}
    (u_n,u_{nx}^2\,dx)\to(u,u_x^2\,dx)\text{ in }\D.
  \end{equation*}
\end{proof}

\begin{proposition}
  \label{prop:convlinf}
  Let $(u_n,\mu_n)$ be a sequence in $\D$ that
  converges to $(u,\mu)$ in $\D$. Then
  \begin{equation*}
    u_n\to u\text{ in }L^\infty(\Real).
  \end{equation*}
\end{proposition}

\begin{proof}
  Let $X_n=(y_n,U_n,H_n)=L(u_n,\mu_n)$ and
  $X=(y,U,H)=L(u,\mu)$, see \eqref{eq:Ldef}. By
  the definition of the metric $d_\D$, we have
  $\lim_{n\to\infty} d(X_n,X)=0$. We claim that
  \begin{equation}
    \label{eq:convLinfty}
    X_n\to X\text{ in }L^\infty(\Real).
  \end{equation}
  The proof of this claim follows the same lines
  as the proof of Lemma \ref{lem:proofdist}. For
  any $\epsi>0$, there exists $N$ such that for
  any $n\geq N$ there exist a path $X^n\in\C_0$
  joining $X_n$ and $X$ such that
  \begin{equation}
    \label{eq:deqepsi1}
    \int_0^1\tnorm{X^n_s}_{X^n(s)}\,ds\leq\frac{\epsi}2.
  \end{equation}
  We have the decomposition
  \begin{equation}
    \label{eq:decXsg1}
    X^n_s(s,\xi)=g^n(s,\xi)X^n_\xi(s,\xi)+R^n(s,\xi).
  \end{equation}
  In the same way that we obtained
  \eqref{eq:decR1R2}, we now obtain
  \begin{equation}
    \label{eq:decR1R21}
    g^n(s,\xi)=R_1^n(s,\xi)+R_3^n(s,\xi).
  \end{equation}
  and it follows that
  \begin{align*}
    \norm{X_n-X}_{L^\infty(\Real)}
    &\leq\int_0^1\norm{X^n_s(s,\cdot)}_{L^\infty}\,ds\\
    &\leq\int_0^1(\norm{g^n(s,\dott)}_{L^\infty}+\norm{R^n(s,\dott)}_{L^\infty})\,ds\\
    &\leq
    2\int_0^1\norm{R^n(s,\dott)}_{L^\infty}\,ds\quad\text{ (by \eqref{eq:decR1R21})}\\
    &\leq
    2\int_0^1\norm{R^n(s,\dott)}_B\,ds=2\int_0^1\tnorm{X^n(s,\dott)}_{X(s)}\,ds
    \leq\epsi.
  \end{align*}
  and this concludes the proof of the claim
  \eqref{eq:convLinfty}. The rest of the proof is
  similar to the proof in \cite[Proposition
  5.2]{HolRay:07}. We reproduce it here for the
  sake of completeness. For any $x\in\Real$, there
  exists $\xi_n$ and $\xi$, which may not be
  unique, such that $x=y_n(\xi_n)$ and
  $x=y(\xi)$. We set $x_n=y_n(\xi)$. We have
  \begin{equation}
    \label{eq:udec}
    u_n(x)-u(x)=u_n(x)-u_n(x_n)+U_n(\xi)-U(\xi)
  \end{equation}
  and
  \begin{align}
    \notag
    \abs{u_n(x)-u_n(x_n)}&=\abs{\int_{\xi}^{\xi_n}U_{n,\xi}(\eta)\,d\eta}\\
    \notag
    &\leq\sqrt{\xi_n-\xi}\left(\int_{\xi}^{\xi_n}U_{n,\xi}^2\,d\eta\right)^{1/2}&\text{(Cauchy--Schwarz)}\\
    \notag
    &=\sqrt{\xi_n-\xi}\left(\int_{\xi}^{\xi_n}y_{n,\xi}H_{n,\xi}\,d\eta\right)^{1/2}&\text{(from
      \eqref{eq:defFGal3})}\\ \notag
    &\leq\sqrt{\xi_n-\xi}\sqrt{\abs{y_n(\xi_n)-y_n(\xi)}}&\text{(since
      $H_{n,\xi}\leq1$)}\\ \notag
    &=\sqrt{\xi_n-\xi}\sqrt{y(\xi)-y_n(\xi)}\\
    \label{eq:undif}
    &\leq \sqrt{\xi_n-\xi}\norm{y-y_n}_{L^\infty(\Real)}^{1/2}.
  \end{align}
  From \eqref{eq:Ldef1}, one can prove that
  \begin{equation*}
    \abs{y(\xi)-\xi}\leq\mu(\Real)
  \end{equation*}
  and it follows that
  \begin{equation*}
    \abs{\xi_n-\xi}\leq
    2\mu_n(\Real)+\abs{y_n(\xi_n)-y_n(\xi)}
    =2 H_n(\infty)+\abs{y(\xi)-y_n(\xi)}
  \end{equation*}
  and, therefore, since $H_n\to H$ and $y_n\to y$ in
  $L^\infty(\Real)$, $\abs{\xi_n-\xi}$ is bounded by a
  constant $C$ independent of $n$. Then,
  \eqref{eq:undif} implies
  \begin{equation}
    \label{eq:undif2}
    \abs{u_n(x)-u_n(x_n)}\leq C
    \norm{y-y_n}_{L^\infty(\Real)}^{1/2}.
  \end{equation}
  Since $y_n\to y$ and $U_n\to U$ in $L^\infty(\Real)$, it
  follows from \eqref{eq:udec} and \eqref{eq:undif2}
  that $u_n\to u$ in $L^\infty(\Real)$.

\end{proof}

\bibliographystyle{plain} 

\end{document}